\documentclass[a4paper,12pt]{article}
\begin{document}

\title{Residues of functions of Cayley-Dickson variables
and Fermat's last theorem.}
\author{Ludkovsky S.V.}
\date{25 August 2010}
\maketitle
\begin{abstract}
Function theory of Cayley-Dickson variables is applied to Fermat's
last theorem. For this the homotopy theorem, Rouch\'e's theorem and
residues of meromorphic functions over Cayley-Dickson algebras are
used. A special meromorphic function of Cayley-Dickson variables is
constructed and its properties are investigated.
\end{abstract}

\section{Introduction.}
\par Analytical methods are frequently used in number theory,
particularly of complex analysis \cite{hasse,irelrosb,manpanchb}. It
is logical to use hypercomplex over Cayley-Dickson algebras analysis
to problems of number theory. One of such interesting objects is
Fermat's last theorem \cite{irelrosb,manpanchb,ribenb}. Its existing
proof is very long and complicated. This article is devoted to a
rather short and clear demonstration of Fermat's last theorem with
the help of (super-)analysis over Cayley-Dickson algebras.
\par Such algebras have a long history, because quaternions were
first introduced by W.R. Hamilton in 1843. He had planned to use
them for problems of mechanics and mathematics \cite{hamilt}. Their
generalization known as the octonion algebra was introduced by J.T.
Graves and A. Cayley in 1843-45. Then Dickson had investigated more
general algebras known now as the Cayley-Dickson algebras
\cite{baez,dickson,kansol}.
\par This work continues previous articles of the
author. In those articles (super)-differentiable functions of
Cayley-Dickson variables and their non-commutative line integrals
were investigated \cite{ludfov,ludoyst,lujmsnfpcd,lufejmsqcf}.
Meromorphic functions of Cayley-Dickson variables, their arguments
and residues were studied in the papers
\cite{ludoyst,ludfov,lufjmsrf}. Super-differentiability or
$z$-differentiability (or $(z,z^*)$- or $z^*$-differentiability)
over Cayley-Dickson algebras ${\cal A}_r$ is a specific
differentiation of the algebra of locally converging formal power
series on an open or canonical closed domain $U$ in ${\cal A}_r$
with prescribed order of multiplication of the Cayley-Dickson
variable $z\in U$ (or $(z, z^*)\in U$ or $z^*\in U$ respectively)
and constants from ${\cal A}_r$ which are ordered in each addendum,
where $z^* = \tilde z$ denotes the conjugated Cayley-Dickson number
$z$. Each locally analytic function is considered on an open
neighborhood of a point of a domain as a given phrase by $z$ (or
$(z, z^*)$ or $z^*$ correspondingly) that to preserve algebraic
features, because functions in the set-theoretic sense do not bear
an algebraic information.
\par The Cayley-Dickson algebras ${\cal A}_r$ have the even
generator $i_0=1$ and the purely imaginary odd generators
$i_1,...,i_{2^r-1}$, $2\le r$, $i_k^2=-1$ and $i_0i_k=i_k$ and
$i_ki_l= - i_li_k$ for each $1\le k\ne l$. For $3\le r$ the
multiplication of these generators is generally non-associative, so
they form not a group, but a non-commutative quasi-group with the
property of alternativity $i_k (i_ki_l) = (i_k^2)i_l$ and
$(i_li_k)i_k=i_l(i_k^2)$ instead of associativity.
\par The purpose of this paper
is in application of developed earlier technique of residues of
meromorphic functions of Cayley-Dickson variables to Fermat's last
theorem. For this the argument principle, homotopy theorem,
Rouch\'e's theorem and residues of meromorphic functions over
Cayley-Dickson algebras are used. A special meromorphic function of
Cayley-Dickson variables is constructed and its properties are
investigated. It is demonstrated that this approach is effective.
\par Main results of this paper are obtained for the first time.
In this article notations and definitions of previous papers cited
above are used.

\section{Special meromorphic function and its poles.}
\par To avoid misunderstanding we first recall some notations
and basic facts in \S \S 1-6 from \cite{ludfov,ludoyst,lufjmsrf}.

\par {\bf 1.} Cayley-Dickson algebras ${\cal A}_r$ form the sequence
so that ${\cal A}_{r+1}$ is obtained from the preceding ${\cal A}_r$
with the help of the so called doubling procedure
\cite{baez,kansol,kurosh}. Therefore, the natural embeddings ${\cal
A}_r \hookrightarrow {\cal A}_{r+1}\hookrightarrow ... $ are
induced. It is convenient to put: ${\cal A}_0 = {\bf R}$ for the
real field, ${\cal A}_1 = {\bf C}$ for the complex field, ${\cal
A}_2 = {\bf H}$ denotes the quaternion skew field, ${\cal A}_3 =
{\bf O}$ is the octonion algebra, ${\cal A}_4$ denotes the sedenion
algebra. The quaternion skew field ${\bf H}$ is associative, but
non-commutative. The octonion algebra ${\bf O}$ is the alternative
division algebra with the multiplicative norm. The sedenion algebra
and Cayley-Dickson algebras of higher order $r\ge 4$ are not
division algebras and have not any non-trivial multiplicative norm.
Nevertheless, they are power-associative, that is $z^nz^m=z^{n+m}$
for any natural numbers $n$ and $m$ and each Cayley-Dickson number
$z$, where $z^n=(z(...(zz)...))$ is the $n$-th power of $z$ (see
also \cite{baez,dickson,kansol}). The norm in the Cayley-Dickson
algebra ${\cal A}_r$ is defined by the equality $|z|^2=zz^*$.

\par We recall the doubling procedure
for the Cayley-Dickson algebra ${\cal A}_{r+1}$ from ${\cal A}_r$,
because it is frequently used. Each Cayley-Dickson number $z\in
{\cal A}_{r+1}$ is written in the form $z=\xi + \eta {\bf l}$, where
${\bf l}^2=-1$, ${\bf l}\notin {\cal A}_r$, $\xi , \eta \in {\cal
A}_r$. The addition of such numbers is componentwise. The conjugate
of any Cayley-Dickson number $z$ is prescribed by the formula:
\par $(1)$ $z^* := {\tilde z} := \xi ^* - \eta {\bf l}$. \\
The multiplication in ${\cal A}_{r+1}$ is defined by the
following equation:
\par $(2)$ $(\xi + \eta {\bf l})(\gamma +\delta {\bf l})=(\xi \gamma
-{\tilde {\delta }}\eta )+(\delta \xi +\eta {\tilde {\gamma }}){\bf l}$ \\
for each $\xi $, $\eta $, $\gamma $, $\delta \in {\cal A}_r$, $z :=
\xi +\eta {\bf l}\in {\cal A}_{r+1}$, $\zeta :=\gamma +\delta {\bf
l} \in {\cal A}_{r+1}$.
\par The basis of ${\cal A}_r$ over $\bf R$ is denoted by
${\bf b}_r := {\bf b} := \{ 1, i_1,...,i_{2^r-1} \} $, where
$i_s^2=-1$ for each $1\le s \le 2^r-1$, $i_{2^r}:={\bf l}$ is the
additional element of the doubling procedure of ${\cal A}_{r+1}$
from the preceding Cayley-Dickson algebra ${\cal A}_r$. Their
enumeration can be chosen as $i_{2^r+m}= i_m{\bf l}$ for each
$m=1,...,2^r - 1$, $i_0:=1$. This implies that $\xi {\bf l} = {\bf
l} \xi ^*$ for each $\xi \in {\cal A}_r$, when $1\le r$.

\par {\bf 2. Remarks and notations.}
\par The family of all ${\cal A}_r$ locally $z$-analytic functions $f(z)$
on a domain $U$ in ${\cal A}_r$ with values in the Cayley-Dickson
algebra ${\cal A}_r$ is denoted by ${\cal H}(U,{\cal A}_r)$ or
$C_z^{\omega }(U,{\cal A}_r)$.
\par To rewrite a function from real variables $z_j$
in the $z$-representation the following identities are used:
\par $(1)$ $z_j=(-zi_j+ i_j(2^r-2)^{-1} \{ -z
+\sum_{k=1}^{2^r-1}i_k(zi_k^*) \} )/2$ \\ for each
$j=1,2,...,2^r-1$, $$(2)\quad z_0=(z+ (2^r-2)^{-1} \{ -z +
\sum_{k=1}^{2^r-1}i_k(zi_k^*) \} )/2,$$  where $2\le r\in \bf N$,
$z$ is a Cayley-Dickson number decomposed as
\par $(3)$ $z=z_0i_0+...+z_{2^r-1}i_{2^r-1}\in {\cal
A}_r$, $z_j\in \bf R$ for each $j$, $i_k^* = {\tilde i}_k = - i_k$
for each $k>0$, $i_0=1$, since $i_k(i_0i_k^*)=i_0=1$, $i_k(i_ji_k^*)
= - i_k(i_k^*i_j) = - (i_ki_k^*) i_j = - i_j$ for each $k\ge 1$ and
$j\ge 1$ with $k\ne j$ (shortly $k\ne j\ge 1$), $i_k (i_ki_k^*)
=i_k$ for each $k\ge 0$.

\par {\bf 3. Notation.} If $f: U\to {\cal A}_r$ is either
$z$-differentiable or $\tilde z$-differentiable at $a\in U$ or on
$U$, then we can write also $D_{\tilde z}$ instead of $\partial
_{\tilde z}=\partial /\partial {\tilde z}$ and $D_z$ instead of
$\partial _z=\partial /\partial z$ at $a\in U$ or on $U$
respectively in situations, when it can not cause a confusion, where
$U$ is an open domain in ${\cal A}_r$.

\par {\bf 4. Proposition.} {\it A function $f: U\to {\cal A}_r$ is
$z$-differentiable at a point $a\in U$ if and only if $F$ is
Fr\'echet differentiable at $a$ and $\partial _{\tilde
z}f(z)|_{z=a}=0$. If $f$ is $z$-super-differentiable on $U$, then
$f$ is $z$-represented on $U$.  A $(z,{\tilde z})$-differentiable
function $f$ at $a\in U$ is $z$-differentiable at $a\in U$ if and
only if $D_{\tilde z}f(z,{\tilde z})|_{z=a}=0$.}
\par {\bf Proof.} For each canonical closed compact set
$U$ in ${\cal A}_r$ the set of all polynomial by $z$ functions is
dense in the space of all continuous on $U$ Fr\'echet differentiable
functions on $Int (U)$ relative to the compact-open topology due the
generalization of Stone-Weierstrass' theorem over the Cayley-Dickson
algebras (see also \cite{ludfov,ludoyst}).
\par As usually a set $A$ having structure of an $\bf R$-linear space
and having distributive multiplications of its elements on
Cayley-Dickson numbers $z\in {\cal A}_v$ from the left and from the
right is called a left- and right-${\cal A}_v$-module (or vector
space over ${\cal A}_v$ if this terminology can not cause any
confusion). \par For two vector spaces $A$ and $B$ over ${\cal A}_v$
one can consider their ordered tensor product $A\otimes B$ over
${\cal A}_v$ consisting of elements $a\otimes b := (a,b)$ such that
$a\in A$ and $b\in B$, $\alpha (a,b)=(\alpha a,b)$ and $(a,b)\beta =
(a,b\beta )$ for each $\alpha , \beta \in {\cal A}_v$, $(a_1\otimes
b_1)(a_2\otimes b_2)= a_1a_2\otimes b_1b_2$ for each $a_1, a_2\in A$
and $b_1, b_2\in B$. In the aforementioned respect $A\otimes B$ is
the $\bf R$-linear space and at the same time left and right module
over ${\cal A}_v$. Then $A\otimes B$ has the structure of the vector
space over ${\cal A}_v$. By induction consider tensor products $\{
C_1\otimes C_2\otimes ... \otimes C_n \} _{q(n)}$, where
$C_1,...,C_n \in \{ A, B \} $, $q(n)$ indicates on the order of
tensor multiplications in the curled brackets $ \{ * \} $. \par For
two ${\cal A}_v$-vector spaces $V$ and $W$ their direct sum $V\oplus
W$ is the ${\cal A}_v$-vector space consisting of all elements
$(a,b)$ with $a\in V$ and $b\in W$ such that $\alpha (a,b)= (\alpha
a,\alpha b)$ and $(a,b)\beta =(a\beta ,b\beta )$ for each $\alpha $
and $\beta \in {\cal A}_v$. Therefore, the direct sum of all
different tensor products $\{ C_1\otimes C_2\otimes ... \otimes C_n
\} _{q(n)}$, which are $\bf R$-linear spaces and left and right
modules over ${\cal A}_v$, provides the minimal tensor space
$T(A,B)$ generated by $A$ and $B$.
\par Operators $\partial _z$ and $\partial _{\tilde z}$ are uniquely defined
on $C^{\omega }_z(U,{\cal A}_r)$ and $C^{\omega }_{\tilde z}
(U,{\cal A}_r)$ in terms of phrases, hence they are unique on the
tensor space $T(C^{\omega }_z(U,{\cal A}_r),C^{\omega }_{\tilde
z}(U,{\cal A}_r))$, which is dense in $C^{\omega }_{z,\tilde
z}(U,{\cal A}_r)$, since $C^{\omega }_{z,\tilde z}(U,{\cal A}_r):=
C^{\omega }_{\mbox{ }_1z,\mbox{ }_2z}(U^2,{\cal A}_r) |_{\mbox{
}_1z=z,\mbox{ }_2z=\tilde z}$. Therefore, operators $\partial _z$
and $\partial _{\tilde z}$ are uniquely defined on $C^{\omega
}_{z,\tilde z}(U,{\cal A}_r)$.
\par If there is a product $fg$ of two phrases $f$ and $g$ from
$C^{\omega }_{z,\tilde z}(U,{\cal A}_r)$, then if it is reduced to a
minimal phrase $\xi $, then it is made with the help of
$z^nz^m=z^{n+m}$ and ${\tilde z}^n{\tilde z}^m={\tilde z}^{n+m}$ and
identities for constants in ${\cal A}_r$, since no any shortening
related with their permutation $z{\tilde z}={\tilde z}z$ or
substitution of $z$ on $\tilde z$ or $\tilde z$ on $z$, for example,
using the identity ${\tilde z}=l(zl^*)$ is not allowed in $C^{\omega
}_{z, \tilde z}(U,{\cal A}_r)$ in accordance with our convention in
\S 2.1 \cite{ludfov,ludoyst}, since $C^{\omega }_{z,\tilde
z}(U,{\cal A}_r):=C^{\omega }_{\mbox{ }_1z, \mbox{ }_2z}(U^2,{\cal
A}_r)|_{\mbox{ }_1z=z,\mbox{ }_2z=\tilde z}$ and in $C^{\omega
}_{\mbox{ }_1z,\mbox{ }_2z}(U^2,{\cal A}_r)$ variables $\mbox{ }_1z$
and $\mbox{ }_2z$ do not commute, $\mbox{ }_1z$ and $\mbox{ }_2z$
are different variables which are not related. Therefore, $\partial
_z\xi .h=(\partial _zf.h)g+f(\partial _zg.h)$ and $\partial _{\tilde
z}\xi .h=(\partial _{\tilde z}f.h)g+ f(\partial _{\tilde z}g.h)$,
hence $\partial _z$ and $\partial _{\tilde z}$ are correctly
defined. \par We consider the super-differentiability, but in
accordance with our convention above for short we simply write
differentiability over ${\cal A}_r$.  \par We can use a $\delta
$-approximation for each $\delta >0$ of $Dg(z).h$ on a sufficiently
small open subset $V$ in $U$ such that $z\in V$ and $|h|\le 1$ by
functions $\zeta _n$ polynomial in $z$ and $\bf R$-homogeneous
${\cal A}_r$-additive in $h$, we also use the partition of unity in
$U$ by $C^{\omega }_z$-functions. Then consider functions $\xi _n$
with ${\xi _n}'.h$ corresponding to $\zeta _n$ for each canonical
closed compact subset $W$ in $U$, since from each open covering of
$W$ we can choose a finite sub-covering of $W$.
\par Suppose that $f$ is $z$-differentiable at a point $a$.
The derivative $f'(z)$ is the ${\bf R}$-linear ${\cal A}_r$-additive
operator, so to it  an $\bf R$-linear operator on the Euclidean
space $\bf R^{2^rn}$ corresponds (see \S 2.1, 2.2
\cite{ludfov,ludoyst}). Then $f(a+h)-f(a)=\partial _af(a).h+
\epsilon (h)|h|$ and $\partial _{\tilde z}f(z)|_{z=a}=0$, since
generally for a $(z,{\tilde z})$-differentiable function
$f(a+h)-f(a)=(\partial _af(a)).h+(\partial _{\tilde a}f(a)).h
+\epsilon (h)|h|$, where $\epsilon (h)$ is continuous by $h$ and
$\epsilon (0)=0$. \par Vice versa, if $F$ is Fr\'echet
differentiable and $\partial _{\tilde z}f(z)|_{z=a}=0$, then
expressing $z_ji_j$ for each $j=0,1,...,2^r-1$ through linear
combinations of $z$ with multiplication on constant coefficients
from ${\cal A}_r$ on the left and on the right in accordance with
Formulas 2$(1-3)$ we get the increment of $f$ as above.
\par The last statement of this proposition follows from
Definition 2.2 \cite{ludfov,ludoyst}.

\par {\bf 5. Proposition.} {\it Let $g: U\to {\cal A}_r$, $r\ge 2$,
and $f: W\to {\cal A}_r$ be two differentiable functions on $U$ and
$W$ respectively such that $g(U)\subset W$, $U$ and $W$ are open in
${\cal A}_r$, where $f$ and $g$ are simultaneously either
$(z,{\tilde z})$, or $z$, or $\tilde z$-differentiable. Then the
composite function $f\circ g(z) := f(g(z))$ is differentiable on $U$
and
\par $(Df\circ g(z)).h = (Df(g)).((Dg(z)).h)$ \\
for each $z\in U$ and each $h\in {\cal A}_r$, and hence $f\circ g$
is of the same type of differentiability as $f$ and $g$.}
\par {\bf Proof.} Theorems 2.11, 2.15, 2.16, 3.10 and
Corollary 2.13 \cite{ludfov,ludoyst} establish the equivalence of
notions of ${\cal A}_r$-holomorphic and ${\cal A}_r$ locally
$z$-analytic classes of functions on open domains in ${\cal A}_r$.
\par In view of these results it is sufficient to prove this Proposition
on open domains $U$ and $W$, where  $g(U)\subset W$ with $U =
g^{-1}(W)$, if others conditions are the same, since
$z$-differentiability is equivalent with the local $z$-analyticity
(that is in the $z$-representation). Indeed, the composition of two
locally $z$-analytic functions $f\circ g$ with such domains is the
locally $z$-analytic function, analogously for
$\tilde z$ and $(z,{\tilde z})$ differentiability.
\par Since $g$ is differentiable, the function $g$ is continuous and $g^{-1}(W)$
is open in ${\cal A}_r$. In view of Proposition 4 above if $f$ and
$g$ are simultaneously either $z$-differentiable or $\tilde
z$-differentiable, then either $\partial _{\tilde z}f=0$ and
$\partial _{\tilde z}g=0$ or $\partial _zf=0$ and $\partial _zg=0$
correspondingly on their domains. \par Consider the increment of the
composite function $$f\circ g(z+h)-f\circ
g(z)=(Df(g))|_{g=g(z)}.(g(z+h)-g(z))+ \epsilon _f(\eta )|\eta |,$$
where $\eta =g(z+h)-g(z)$, $g(z+h)-g(z)=(Dg(z)).h+\epsilon _g(h)|h|$
(see \S 2.2 \cite{ludfov,ludoyst}). Since the derivative $Df$ is
${\cal A}_r$-additive and $\bf R$-homogeneous
(and continuous) operator on ${\cal A}_r$, we have \\
$f\circ g(z+h)-f\circ g(z)=(Df(g))|_{g=g(z)}.((Dg(z)).h) +\epsilon
_{f\circ g}(h)|h|$, where
$$\epsilon _{f\circ g}(h)|h|:=\epsilon _f((Dg(z)).h+\epsilon _g(h)|h|)
|(Dg(z)).h+\epsilon _g(h)|h|| $$
$$+[(Df(g))|_{g=g(z)}.(\epsilon _g(h))]|h|),$$
$$|(Dg(z)).h+\epsilon _g(h)|h||\le
[\| Dg(z) \| + |\epsilon _g(h)|] |h|\mbox{, hence}$$
$$|\epsilon _{f\circ g}(h)|\le |\epsilon _f((Dg(z)).h+\epsilon _g(h)|h|)|
[\| Dg(z) \| + |\epsilon _g(h)|] + \| (Df(g))|_{g=g(z)} \| |\epsilon
_g(h)|$$ and inevitably $\lim_{h\to 0} \epsilon _{f\circ g}(h)=0$.
Moreover, $\epsilon _{f\circ g}(h)$ is continuous in $h$, since
$\epsilon _g$ and $\epsilon _f$ are continuous functions, $Df$ and
$Dg$ are continuous operators. Evidently, if $\partial _{\tilde
z}f=0$ and $\partial _{\tilde z}g=0$ on domains of $f$ and $g$
respectively, then $\partial _{\tilde z}f\circ g=0$ on $V$, since
$D=\partial _z+\partial _{\tilde z}$.

\par Suppose now that there are phrases corresponding to $f_n$ and $g_n$
denoted by $\mu _n$ and $\nu _n$ such that $f_n$ and $g_n$ uniformly
converge to $f$ and $g$ respectively on each bounded canonical
closed subset in $W$ and $U$ from a family $\cal W$ or $\cal U$
respectively, where $\cal W$ and $\cal U$ are coverings of $W$ and
$U$ correspondingly. While $D_z\mu _n$ and $D_z\nu _n$ are
fundamental sequences uniformly on each bounded canonical closed
subset $P\in \cal W$ and $Q\in \cal U$ correspondingly relative to
the operator norm (see Definitions in \S  2 \cite{ludfov,ludoyst}).
Then the sequence $\mu _n\circ \nu _n$ converges on each bounded
canonical closed subset $Q_1$ such that $Q_1\subset Q\cap
g^{-1}(P)$, where $P\in \cal W$ and $Q\in \cal U$ in $U$, when $n$
tends to the infinity. Moreover, $D_y\mu _n(y,{\tilde y}).(D_z\nu
_n(z,{\tilde z}))|_{y=\nu _n(z,{\tilde z})}$ and $D_{\tilde y}\mu
_n(y,{\tilde y}).(D_z{\tilde {\nu }} _n(z,{\tilde z}))|_{y=\nu
_n(z,{\tilde z})}$ are  the fundamental sequences of operators on
each bounded canonical closed subset $Q_1$ so that $Q_1\subset Q\cap
g^{-1}(P)$, where $P\in \cal W$ and $Q\in \cal U$. The family $
{\cal Q} := \{ Q_1 \} $ specified above evidently is the covering of
$V$.
\par It remains to verify, that $D_z[\mu _n\circ \nu _n(z,{\tilde z})] =
[(D_y\mu _n(y,{\tilde y}).(D_z\nu _n(z,{\tilde z})) + (D_{\tilde
y}\mu _n(y,{\tilde y}).(D_z{\tilde {\nu }} _n(z,{\tilde
z}))]|_{y=\nu _n(z,{\tilde z})}$.
\par Since the derivation operator $D_z$ by $z$ is
$\bf R$-homogeneous and ${\cal A}_r$-additive,
it is sufficient to verify this in the case $D_z[\eta \circ \psi
(z,{\tilde z})]$ locally in balls, where both series uniformly
converge. Here the phrase is written as:
$$(1)\quad \eta = \eta (z, {\tilde z}) =
\sum_k \{ A_k, z, {\tilde z } \}_{q(k)} := \{
a_{k,1}\overbrace{z^{k_1}_{l_1}}...a_{k,p}\overbrace{z^{k_p}_{l_p}}
\}_{q(k)} ,$$ $k=(k_1,...,k_p)$, $p=p(k)\in \bf N$, $0\le k_j\in \bf
Z$ and $l_j \in \{ 1, 2 \} $ for each $j$, $a_{k,1},...,a_{k,p}\in
{\cal A}_r$ are constants, $\overbrace{z_1}=z$,
$\overbrace{z_2}=\tilde z$, $z$ is the Cayley-Dickson variable, a
vector $q$ indicates on an order of the multiplication, $A_k =
[a_{k,1},...,a_{k,p}]$, also
$$(2)\quad \psi = \psi (z, {\tilde z}) =
\sum_m \{ B_m, z, {\tilde z} \}_{q(m)} ,$$ where $m=(m_1,...,m_s)$,
$B_m = [b_{m,1},...,b_{m,s}]$, $b_{m,j}\in {\cal A}_r$ is a constant
for each $m, j$, also $s\in \bf N$, $0\le m_j\in \bf N$ for each
$j$. We have the identities
$$(3)\quad D_z\overbrace{z_l} = \delta _{l,1} {\bf 1},\quad D_{\tilde z}
\overbrace{z_l} = \delta _{l,2} {\tilde {\bf 1}},$$ where ${\bf 1}h
:= h$ and ${\tilde {\bf 1}}h := \tilde h$ for each $h\in {\cal
A}_r$. Using shifts $z\mapsto z-\zeta $ we can consider that the
series are decomposed around a zero point. If a series $\eta $ is
uniformly converging on a canonical closed subset $Q_1$ as above,
then
$$(4)\quad D_z\eta =\sum_k D_z \{ a_{k,1}\overbrace{z
^{k_1}_{l_1}}...a_{k,p}\overbrace{z^{k_p}_{l_p}} \}_{q(k)} $$ is
also uniformly converging on $Q_1$ in accordance with our
supposition about convergence of phrases and their derivatives.
\par In suitable $Q_1$ we deduce from Formula
$(4)$ that
$$(5)\quad D_z[\eta \circ \psi (z,{\tilde z})] = D_z
\sum_k \{ a_{k,1}\overbrace{\psi
^{k_1}_{l_1}}...a_{k,p}\overbrace{\psi ^{k_p}_{l_p}} \}_{q(k)} $$ $$
= \sum_k D_z \{ a_{k,1}\overbrace{\psi
^{k_1}_{l_1}}...a_{k,p}\overbrace{\psi ^{k_p}_{l_p}} \}_{q(k)}$$
$$ = \sum_k [ \{ a_{k,1}(D_z\overbrace{\psi
^{k_1}_{l_1}})...a_{k,p}\overbrace{\psi ^{k_p}_{l_p}} \}_{q(k)}+...+
\{ a_{k,1}\overbrace{\psi ^{k_1}_{l_1}}...a_{k,p}(D_z\overbrace{\psi
^{k_p}_{l_p}}) \}_{q(k)} ]$$  $$ = \sum_k \sum_{\mbox{
}_1m,....,\mbox{ }_pm; p=p(k)} [ \{ a_{k,1}(D_z(\overbrace{ \{
B_{\mbox{ }_1m}, z, {\tilde z} \} _{q(\mbox{ }_1m)} } )
^{k_1}_{l_1})...a_{k,p} (\overbrace{ \{ B_{\mbox{ }_pm}, z, {\tilde
z} \} _{q(\mbox{ }_pm)} } ) ^{k_p}_{l_p} \}_{q(k)}$$  $$+...+ \{
a_{k,1} (\overbrace{ \{ B_{\mbox{ }_1m}, z, {\tilde z} \} _{q(\mbox{
}_1m)} } ) ^{k_1}_{l_1}...a_{k,p}(D_z(\overbrace{ \{ B_{\mbox{
}_pm}, z, {\tilde z} \} _{q(\mbox{ }_pm)} } )^{k_p}_{l_p}) \}_{q(k)}
.$$  On the other hand, we have $$D_z \{ B_m, z, {\tilde z }
\}_{q(m)} = \{ b_{m,1}(D_z\overbrace{z
^{m_1}_{l_1}})...b_{m,s}\overbrace{z^{k_s}_{l_s}} \}_{q(m)} + ... +
\{ b_{m,1}\overbrace{z
^{m_1}_{l_1}}...b_{m,s}(D_z\overbrace{z^{k_s}_{l_s}}) \}_{q(k)}$$
due to the Leibnitz rule, hence
$$ D_z\overbrace{\psi ^p} = \sum_{\mbox{
}_1m,....,\mbox{ }_pm} [(... (D_z\overbrace{ { \{ B_{\mbox{ }_1m},
z, {\tilde z} \}  }_{q(\mbox{ }_1m)} })...)\overbrace{ { \{
B_{\mbox{ }_pm}, z, {\tilde z} \} } _{q(\mbox{ }_pm)} } $$
$$+...+ (...( \overbrace{ { \{ B_{\mbox{ }_1m}, z, {\tilde z} \} }
_{q(\mbox{ }_1m)} } \overbrace{ { \{ B_{\mbox{ }_2m}, z, {\tilde z}
\} } _{q(\mbox{ }_2m)} } )...)(D_z\overbrace{ { \{ B_{\mbox{ }_pm},
z, {\tilde z} \} } _{q(\mbox{ }_pm)} } ) ]$$
$$ = (D_y\overbrace{y^p}).(D_z\psi ) + (D_{\tilde y}\overbrace{y^p}).(D_z
{\tilde {\psi }})$$ due to Formulas $(3)$, where $y=\psi (z,{\tilde
z})$. Thus Formulas $(5,6)$ imply that
$$(7)\quad D_z (\eta \circ \psi (z,{\tilde z})) = [ (D_y\eta
(y,{\tilde y})).(D_z\psi (z, {\tilde z})) + (D_{\tilde y}\eta
(y,{\tilde y})).(D_z{\tilde {\psi }}(z, {\tilde z})) ]|_{y=\psi (z,
{\tilde z})}.$$ Quite analogously or using the conjugation one
deduces that
$$(8)\quad D_{\tilde z} (\eta
\circ \psi (z,{\tilde z})) = [ (D_{\tilde y}\eta (y,{\tilde
y})).(D_{\tilde z} \psi (z, {\tilde z})) + (D_y\eta (y,{\tilde
y})).(D_{\tilde z}{\tilde {\psi }}(z, {\tilde z})) ]|_{y=\psi (z,
{\tilde z})}.$$ Particularly, we get
$$(9)\quad D_z (\eta \circ \psi (z)) = (D_y\eta
(y)).(D_z\psi (z))|_{y=\psi (z)},\mbox{ when } D_{\tilde z}\eta
(z,{\tilde z})=0\mbox{ and } D_{\tilde z}\psi (z,{\tilde z})=0$$ and
$$(10)\quad D_{\tilde z} (\eta
\circ \psi ({\tilde z})) = (D_{\tilde y}\eta ({\tilde
y})).(D_{\tilde z}\psi ({\tilde z})) |_{y=\psi ({\tilde z})},\mbox{
when } D_z\eta (z,{\tilde z})=0\mbox{ and } D_z\psi (z,{\tilde
z})=0.$$ Combining Formulas $(7,8)$ and taking into account the
proof given above and applying Proposition 2.3 \cite{ludfov,ludoyst}
we infer the chain rule, when both $\eta $ and $\psi $ are either
the $z$ or $\tilde z$ or $(z,{\tilde z})$-differentiable, that is in
all three considered cases.

\par {\bf 6. Some elementary functions and their non-commutative Riemann surfaces.}
\par In this article some elementary facts about analytic functions
$z^n$, $z^{1/n}$, $\exp (z)$ and $Ln (z)$  of the Cayley-Dickson
variables are used. They were considered in details in previous
works \cite{ludoyst,ludfov,lujmsnfpcd}. Recall that the exponential
function is defined by the power series $$\exp (z) := 1+
\sum_{n=1}^{\infty } z^n/n!$$ converging on the entire
Cayley-Dickson algebra ${\cal A}_r$. It has the periodicity property
$\exp (M(\phi + 2\pi k))=\exp (M\phi )$ for each purely imaginary
Cayley-Dickson number $M$ of the unit norm $|M|=1$ for any real
number $\phi \in {\bf R}$ and each integer number $k\in {\bf Z}$.
The restriction of such exponential function on each complex plane
${\bf C}_M := {\bf R}\oplus M{\bf R}$ coincides with the traditional
complex exponential function. Since the inverse function $Ln (z)$ of
$z=\exp (x)$ is defined on every complex plane ${\bf C}_M\setminus
\{ 0 \} $ with the pricked zero point, the logarithmic function $Ln
(z)$ is defined on ${\cal A}_r\setminus \{ 0 \} $. This logarithmic
function is certainly multi-valued. Consider the bunch of complex
planes ${\bf C}_M$ intersecting by the real line ${\bf R}i_0$ as the
geometric realization of ${\cal A}_r$. Certainly ${\bf C}_M = {\bf
C}_{-M}$, so we take the set \par ${\bf S}_r^+ := \{ M\in {\cal
A}_r:$ $|M|=1,$ $M= M_1 i_1 +...+ M_{2^r-1} i_{2^r-1} ,$
$\mbox{either}$ $M_1>0,$ $\mbox{or}$ $M_1=0$ $\mbox{and}$ $M_2>0$,
$\mbox{or}$...$,\mbox{ or}$ $M_1=...=M_{2^r-2}=0$ $\mbox{and}$
$M_{2^r-1}>0 \} $, where $M_1,...,M_{2^r-1}\in {\bf R}$, $2\le r$.
Then ${\cal A}_r = \bigcup_{M\in {\bf S}_r^+} {\bf C}_M$.
\par If $A$ and $B$ are two subsets in a complete uniform space $X$
and $\theta : A\to B$ is a continuous bijective mapping, the
equivalence relation $a\Upsilon b$ by definition means $b=\theta
(a)$ for $a\in A$ and $b\in B$; or $b\Upsilon a$ means $a=\theta
^{-1}(b)$. When $\theta $ is uniformly continuous, $\theta $ has a
uniformly continuous extension $\theta : cl(A)\to cl (B)$, where
$cl(A)$ denotes the closure of the set $A$ in $X$ (see Theorem
8.3.10 in \cite{eng}). Certainly, the mapping $\theta $ can be
specified by its graph $ \{ (x_1,x_2): ~ x_2 = \theta (x_1), ~ x_1
\in A \} $. We say, that $A$ and $B$ are glued (by $\theta $), if
$B=\theta (A)$ and the natural quotient mapping $\pi : X\to
X/\Upsilon $ is given, where $X/\Upsilon $ denotes the quotient
space (see \S 2.4 \cite{eng}). For the uniformly continuous $\theta
$ this means that the gluing is extended from $A$ onto $cl (A)$.
\par In the complex case to construct the Riemann surface of the
logarithmic function one takes traditionally the complex plane ${\bf
C}$ cut by the set $ Q_1 := \{ z=x+iy\in {\bf C}: ~ x< 0 \} $ and
marking two respective points $x_1$ and $x_2$  of two edges
$Q_{1,1}$ and $Q_{1,2}$ of the cut $Q_1$ arising from each given
point $x<0$, where $i=i_1$. Then one embeds ${\bf C}$ into either
${\bf C}\times {\bf R}$ or ${\bf C}\times {\bf R}i$ and bents the
obtained surface slightly along the perpendicular axis $e_3$ to
${\bf C}$ by neighborhoods of two edges of the cut $Q_1$ and gets
the new surface ${\cal C}$. Taking the countable infinite family of
such surfaces ${\cal C}^j$ with the edges of cuts $Q^j_{1,1}$ and
$Q^j_{1,2}$ and gluing by respective points of the cuts $Q^j_{1,2}$
with $Q^{j+1}_{1,1}$ for each $j$ one gets the Riemann surface of
the logarithmic function, where $j\in {\bf Z}$ (see, for example,
\cite{lavrsch}).
\par Analogous procedure to construct the Riemann surface
is in the cases $r\ge 2$: one cuts ${\cal A}_r$ by $Q_r$ and gets
two edges $Q_{r,1}$ and $Q_{r,2}$ of the cut. This is described
below.
\par If $K$ and $M=KL$ are two purely imaginary Cayley-Dickson
numbers with $|K|=|M|=|L|=1$ so that they are orthogonal $K\perp M$,
that is $Re(KM)=0$, then $(K{\bf R})\oplus (M{\bf R}) = K({\bf
R}\oplus L{\bf R})$ and $L$ is also purely imaginary. Consider any
path $\gamma : [0,1]\to K{\bf R}\oplus M{\bf R}$ winding one time
around zero such that $\gamma (t)\ne 0$ for each $t$. Each $z\in
K{\bf R}\oplus M{\bf R}$ can be written in the polar form $z=|z| K
e^{L\phi } = |z| \exp ( \pi Ke^{L\phi }/2)$, where $\phi = \phi (z)
\in {\bf R}$, since $Ke^{L\phi } = K\cos (\phi ) + (KL) \sin (\phi
)$ due to Euler's formula, hence $|Ke^{L\phi (\gamma (t))}| =1$ for
each $t$ (see Section 3 in \cite{ludoyst,ludfov}). In particular,
$\gamma (t) = |\gamma (t)| \exp (\pi Ke^{L\phi (t)}/2)$. But
$Ke^{L\phi (\gamma (0))} = Ke^{L\phi (\gamma (1))}$ such that the
logarithm $Ln ~ \gamma (t) = Ln ~ |\gamma (t)| + \pi Ke^{L\phi
(\gamma (t))}/2$ does not change its branch, when the path $\gamma
(t)\in (K{\bf R})\oplus (M{\bf R})\setminus \{ 0 \} $ winds around
zero, since $|\pi Ke^{L\phi } /2|=\pi /2 <\pi $. Due to the homotopy
theorem (see \cite{ludoyst,ludfov}) this means that the logarithm
$Ln ~ \gamma (t)$ does not change its branch, when $Re (\gamma
(t))=0$ for each $t$ for the path $\gamma $ winding around zero with
$|\gamma (t)|>0$ for each $t$.
\par The first simple construction for $r\ge 2$ is the following. Take the set
$Q_r := (-\infty ,0){\bf S}_r^+ := \{ z=tx: ~ t\in (-\infty ,0),
x\in {\bf S}_r^+ \} $ and cut ${\cal A}_r$ by $Q_r$. The set $Q_r$
is the union $Q_r = \bigcup_{j=1}^{2^r-1} \Omega _{j,r}$ of subsets
$\Omega _{j,r} := \{ z\in Q_r: ~ z_0=0,...,z_j=0 \} $ so that
$\Omega _{j,r}$ is contained in the boundary of the preceding set
$\Omega _{j,r} \subset \partial \Omega _{j-1,r}$ for each
$j=1,...,2^r-1$, $dim ~ \Omega _{j,r} = dim ~ \Omega _{j-1,r} -1$,
moreover, ${\bf R}{\bf S}_r^+ = {\cal I}_r := \{ z\in {\cal A}_r: ~
Re (z)=0 \} $. Therefore, from each point $z\in Q_r$ two and only
two different points $z_1$ and $z_2$ arise while cutting of ${\cal
A}_r$ by $Q_r$.
\par  It is useful to embed ${\cal A}_r$ either into ${\cal
A}_r\times {\bf R}^{2^r-1}$ or into ${\cal A}_r\times {\cal I}_r$.
Then one marks all pairs of respective points $z_1$ and $z_2$
arising from $z\in Q_r$ after cutting, slightly bents the cut copy
of ${\cal A}_r\setminus \{ 0 \} $ by $(2^r-1)$ axes perpendicular to
${\cal A}_r$ by two neighborhoods of two edges of the cut $Q_r$.
Thus one gets the $2^r$ dimensional surface ${\cal C}_r$ with two
edges $Q_{r,1}$ and $Q_{r,2}$ of the cut. Taking the countable
infinite family of such surfaces ${\cal C}_r^j$ with edges of the
cuts $Q^j_{r,1}$ and $Q^j_{r,2}$, $j\in {\bf Z}$, and gluing
respective points of edges $Q^j_{r,2}$ with $Q^{j+1}_{r,1}$ for each
$j$ one gets the Riemann surface ${\cal R}_r = {\cal R}_{r,Ln}$ of
the logarithmic function $Ln (z): {\cal A}_r\setminus \{ 0 \} \to
{\cal R}_r$. Thus the latter mapping is already univalent with the
image in the Riemann surface (see in details
\cite{ludoyst,ludfov,lujmsnfpcd}).
\par For convenience we attach numbers $1$ and $2$ to faces
in such manner that the winding around zero in the complex plane
${\bf C}_M$ embedded into ${\cal R}_r$ counterclockwise means the
transition through the cut from $Q^j_{r,2}$ to $Q^j_{r,1}$ for each
$M$ in the connected set ${\bf S}_r^+$.
\par For the function $z^{1/n}$ with $n\in {\bf N}$ its Riemann
surface ${\cal R}_{r,z^{1/n}}$ is obtained from $n$ copies of
surfaces ${\cal C}_r^j$, $j=1,...,n$, by gluing the corresponding
points of edges $Q^j_{r,2}$ with $Q^{j+1}_{r,1}$ for $j=1,...,n-1$
and of $Q^n_{r,2}$ with $Q^1_{r,1}$ (see also
\cite{lavrsch,ludoyst,ludfov,lujmsnfpcd}).
\par Another more complicated construction is described below.
Now we take the set \par $Q_r := \bigcup_{j=1}^{2^r-1} P_j$, where
\par $P_j := \{ z\in {\cal A}_r:$ $z=z_0i_0+...+z_{2^r-1}i_{2^r-1}$;
$z_0<0$ $\mbox{and}$ $z_j=0 \} $, where $z_0,...,z_{2^r-1}\in {\bf
R}$, $2\le r$.
\par Let $z=z_0+z'$ be the
Cayley-Dickson number with the negative real part $z_0<0$ and the
imaginary part $Im (z)=z'$, which can be written in the form $z' =
|z'| \exp (\pi Ke^{L\phi (z')}/2 )$ and $z=|z|e^{P\psi }$, where
$K$, $L$ and $P$ are purely imaginary Cayley-Dickson numbers of the
unit norm, $\phi $ and $\psi \in {\bf R}$ are reals, $Re (KL^*)=0$,
$~|KL|=1$. This gives the relation $\cos (\psi ) =z_0/|z|$ so that
the parameter $\psi $ is in the interval $\pi /2 +2\pi k <\psi <3\pi
/2 +2\pi k$ for some integer number $k$. This means that for a
continuous path $\gamma $ contained in the set $Q_r$ the parameter
$\psi (\gamma (t))$ is the continuous function of the real variable
$t\in {\bf R}$ and remains in the same interval $(\pi /2 +2\pi k ,
3\pi /2 +2\pi k)$, consequently, the logarithmic function $Ln \gamma
(t)$ preserves its branch along such path $\gamma (t)$. This shows
that after the first cut along $Q_r$ the obtained sets $Q_{r,1}$ and
$Q_{r,2}$ need not be further cut. Thus the described reason
simplifies the construction of the Riemann surface.
\par Each continuous path $\gamma : [0,1]\to {\cal A}_r$ can be
decomposed as the point-wise sum and as the composition (join) up to
the homotopy satisfying the conditions of the homotopy theorem
\cite{ludoyst,ludfov} of paths $\gamma _{k,l}$ in the planes $({\bf
R}i_k)\oplus ({\bf R}i_l)$ for each $k<l\in \Lambda $ for the
corresponding subset $\Lambda \subset \{ 0, 1,..., 2^r-1 \} $. If
$|\gamma (t)|>0$ for each $t$ we take $\gamma _{k,l}$ with $|\gamma
_{k,l}(t)|>0$ on $[0,1]$ for all $k<l\in \Lambda $. When $\gamma
[0,1]$ does not intersect $Q_r$ one can choose $\gamma _{k,l}$ with
images $\gamma _{k,l}[0,1]$ also non-intersecting with $Q_r$ for all
$k<l\in \Lambda $. Due to the homotopy theorem the logarithm $Ln ~
\gamma (t)$ does not change its branch along such continuous path
$\gamma (t)$, since this is the case for $Ln ~ \gamma _{k,l}(t)$ for
all $k<l\in \Lambda $.
\par For each $z\in P_j\setminus \bigcup_{m, m\ne j} P_m$ cutting by
$Q_r$ gives two points. If $z\in (P_k\cap P_j)\setminus \bigcup_{m,
m\ne k, m\ne j} P_m $ with $k<j$ cutting gives four points which can
be organized into respective pairs after cutting of $P_k$ and then
of $P_j$. This procedure gives pairs $(z_{1,1}; z_{2,2})$ and
$(z_{1,2}; z_{2,1})$. For $z\in (P_{k_1}\cap ... \cap
P_{k_l})\setminus \bigcup_{m, m\ne k_1,..., m\ne k_l} P_m $ consider
pairs of points appearing from the preceding point after cutting of
$P_{k_j}$, $j=1,...,l$ by induction, where $k_1<...<k_l$. One can do
it by induction by all $l=1,...,2^r-1$ and all possible subsets
$1\le k_1<k_2<...<k_l\le 2^r-1$. It produces points denoted by
$z_{a_1,...,a_l}$. Here $a_j=1$ corresponds to the face of the cut
indicated by the condition $z_j\le 0$, while $a_j=2$ corresponds to
the face of the cut with $z_j\ge 0$. Two points $z_{a_1,...,a_l}$
and $z_{b_1,...,b_l}$ form the pair of respective points, when
$a_j+b_j=3$ for each $j$, where $l\in \{ 1,...,2^r-1 \} $ is the
index of such points. To each point $z\in (P_{k_1}\cap ... \cap
P_{k_l})\setminus \bigcup_{m, m\ne k_1,..., m\ne k_l} P_m $ of index
$l\ge 1$ the point $M\in {\bf S}_r^+$ corresponds such that
$M=M_1i_1+...+M_{2^r-1}i_{2^r-1}$ and $z\in {\bf C}_M$, where
$M_{k_1}=0$,...,$M_{k_l}=0$.
\par Then as above after cutting of ${\cal A}_r$ by $Q_r$
one gets the $2^r$ dimensional surface ${\cal C}_r$ with two edges
$Q_{r,1}$ and $Q_{r,2}$ of the cut. For the countable infinite
family of the surfaces ${\cal C}_r^j$ with edges of the cuts
$Q^j_{r,1}$ and $Q^j_{r,2}$, $j\in {\bf Z}$, one glues respective
points of edges $Q^j_{r,2}$ with $Q^{j+1}_{r,1}$ for each $j$. Thus
one gets the Riemann surface ${\cal R}_r = {\cal R}_{r,Ln}$ of the
logarithmic function $Ln (z): {\cal A}_r\setminus \{ 0 \} \to {\cal
R}_r$. Therefore, the latter mapping is univalent with the image in
the Riemann surface (see in details
\cite{ludoyst,ludfov,lujmsnfpcd}). For simplicity one can glue at
first pairs of respective points of index $l=1$ and extend gluing by
continuity on points of index $l>1$ (see above). The first
construction for $r\ge 2$ operates with points of index one only.
That is why it is simpler than the second procedure.
\par The reader can lightly see that geometrically ${\cal R}_{r,f}$
(for both types described above) is the bunch $ \bigcup \{ {\cal
R}_{1,f,M}: M\in {\bf S}_r^+ \} $ of complex Riemann surfaces ${\cal
R}_{1,f,M}$ for the restrictions $f|_{{\bf C}_M}$ on the complex
planes ${\bf C}_M$ of the function $f$, where either $f=Ln$ or
$f=z^{1/n}$ among those functions considered here.
\par The main problem of the multi-dimensional geometry is in
depicting its objects on the two-dimensional sheet of paper so one
uses either projections or sections of the multi-dimensional
surfaces. The bunch interpretation $\bigcup_{M\in {\bf S}_r^+} {\cal
R}_{1,f,M}$ helps to imagine how the $2^r$-dimensional Riemann
surface ${\cal R}_{r,f}$ is organized for $f=Ln$ and $f=z^{1/n}$
with $n\in {\bf N}$.

\par {\bf 7.} In this section a special meromorphic function
of Cayley-Dickson variables is constructed and its properties are
used to demonstrate the following theorem.
\par {\bf Fermat's last theorem.} {\it No three positive
integers $a$, $b$ and $c$ can satisfy the equation $a^n+b^n=c^n$ for
any integer value of $n$ greater than two.}
\par {\bf Proof.} We take the Cayley-Dickson algebra
${\cal A}_v$ with a sufficiently large index $v\ge 3$ (see below)
and choose three different ${\bf R}$-linear embeddings of the
quaternion skew field ${\bf H}={\cal A}_2$ into it. For this we take
the standard basis of generators $\{ i_0,...,i_{2^v-1} \} $ of
${\cal A}_v$. These embeddings use the doubling generators $\{ i_1,
i_2 \} $ for the first; $\{ i_4, i_8 \} $ for the second and $\{
i_{16}, i_{32} \} $ for the third embedding correspondingly and
their different products, since these purely imaginary pairs $\{ M,
K \} $ satisfy the following conditions $Re (MK)=0$ and $MK= - KM$
and $|MK|=|M||K|$, $|M|=|K|=1$ (see also
\cite{baez,dickson,kansol}).
\par The first embedding $\theta _1: {\bf H}\to {\cal A}_r$ is
the usual identical mapping $\theta _1(i_k)=i_k$ for $k=0, 1, 2, 3$.
We certainly have $\theta _l (x)=x$ for each real number $x$ and
$l=1, 2, 3$. The second embedding $\theta _2$ corresponds to ${\bf
H}$ written as $\bigoplus_{k=0}^3 {\bf R} i_{4k}$, where we use
enumeration $i_{2^r+p}=i_pi_{2^r}$ for $r\ge 1$ and $p=1,...,2^r-1$.
The third embedding $\theta _3$ corresponds to ${\bf H}$ with the
basis of generators $\{ i_{16k}: k=0,...,3 \} $, i.e. $\theta
_3({\bf H}) = \bigoplus_{k=0}^3 {\bf R} i_{16k}$.
\par In three copies ${\bf H}_l := \theta _l({\bf H})$, $l=1, 2, 3$
of the quaternion skew field ${\bf H}$ embedded into the
Cayley-Dickson algebra ${\cal A}_v$ we consider contained in them
copies of the complex field ${\bf C}_l = \theta _l({\bf C}) =
\bigoplus_{k=0}^1 {\bf R} i_{2^{2(l-1)}k}$. Then we take the
doubling generators $i_{64}$ and $i_{128}$ and $i_{256}$ such that
$({\bf H}_1+{\bf H}_2 +{\bf H}_3 + {\bf H}_2{\bf H}_3)$ and ${\bf
H}_1 i_{64}$ and ${\bf H}_2 i_{128}$ and ${\bf H}_3 i_{256}$ do not
intersect pairwise, where $({\bf H}_1+{\bf H}_2 +{\bf H}_3 + {\bf
H}_2{\bf H}_3) = \{ z\in {\cal A}_v: ~ z= \mbox{}_1z + \mbox{}_2z +
\mbox{}_3z+ \mbox{}_4z, ~ \mbox{}_lz \in {\bf H}_l, ~ l=1, 2 ,3, ~
\mbox{}_4z\in {\bf H}_2{\bf H}_3 \} $, $~ AB := \{ z\in {\cal A}_v:
~ z= xy, ~ x\in A, ~ y\in B \} $ for subsets $A$ and $B\subset {\cal
A}_v$. Variables in them we denote by $w\in {\bf C}_1$, $x\in {\bf
C}_2$ and $y\in {\bf C}_3$. That is $wi_{64}+xi_{128}+yi_{256}\in
{\cal A}_v$ for each $w, x, y$. Therefore, it is sufficient to take
$v=9$.
\par Two Cayley-Dickson numbers $p$ and $s$ are orthogonal $p\perp
s$ by the definition if an only if $Re (ps^*)=0$, because \par $(1)$
$Re (ps^*)=\sum_{j=0}^{2^v-1} p_js_j$. \\ Next we take $\bf
R$-linear projection mappings. For each $j=0,...,2^v-1$ the $\bf
R$-linear projection operator $\pi _j : {\cal A}_v\to {\bf R}i_j$
exists due to Formulas 2$(1-3)$ so that $\pi _j(z)=i_jz_j=z_ji_j$:
\par $(2)$ $\pi _j(z) = (- i_j (zi_j) - (2^v-2)^{-1} \{ -z
+\sum_{k=1}^{2^v-1}i_k(zi_k^*) \} )/2$
\\ for each $j=1,2,...,2^v-1$, $$(3)\quad \pi _0(z) = (z+
(2^v-2)^{-1} \{ -z + \sum_{k=1}^{2^v-1}i_k(zi_k^*) \} )/2,$$  where
$2\le v\in \bf N$. Here we take $v=9$.  \par Combining into suitable
sums these projection operators one gets ${\bf R}$-linear projection
operators $\upsilon _l: {\cal A}_9\to {\bf C}_l i_{\kappa (l)}$,
where $\kappa (1)=64$, $\kappa (2)=128$ and $\kappa (3) = 256$, $~
\upsilon _1= \sum_{j=0}^1\pi _{j+64}$, $ ~\upsilon _2 =
\sum_{j=0}^1\pi _{4j+128}$, $ ~\upsilon _3 = \sum_{j=0}^1\pi
_{16j+256}$. The latter projection operators induce ${\bf R}$-linear
operators $\tau _l : {\cal A}_{v}\to {\bf C}_l$ so that $\tau
_l(z)=\upsilon _l(z)i_{\kappa (l)}^*$ for each $z\in {\cal A}_{v}$
and $l=1, 2, 3$, particularly, $\tau _1(wi_{\kappa (1)}+z) = w$ for
each $z\perp wi_{\kappa (1)}$, $ ~ \tau _2(z+xi_{\kappa (2)})=x$ for
each $z\perp xi_{\kappa (2)}$ and $\tau _3(z+yi_{\kappa (3)})=y$ for
each $z\perp yi_{\kappa (3)}$.
\par In accordance with Formulas 2$(1-3)$ and $(2,3)$ each mapping
$\tau _l(z)$ has the finite phrase expression of the type \par $(4)$
$\tau _l(z) = \sum_m (\alpha _{l,m}(z\beta _{l,m}))i_{\kappa (l)}^*
$ \\ in the $z$-representation, which we fix, where $z\in {\cal
A}_{v}$, $ ~ \alpha _{l,m}$ and $\beta _{l,m}\in {\cal A}_{v}$ are
Cayley-Dickson constants.
\par With the help of these ${\bf R}$-linear mappings $\tau _l$ we define
the function
\par $(5)$ $g(z) := \exp [ \tau _1(z) \exp (2\pi
L_1 \tau _1^{1/n}(z) ) ]$\\ $ - (\exp [ \tau _2(z) \exp (2\pi L_2
\tau _2^{1/n}(z)) ]) (\exp [ \tau _3(z) \exp (2\pi L_3 \tau
_3^{1/n}(z)) ] ) + | z - \sum_{l=1}^3 \upsilon _l(z) |^2 L_4$
\\ on the Cayley-Dickson algebra ${\cal A}_{v}$, where $L_1=i_2$, $  ~
L_2 = i_{8}$ and $L_3 = i_{32}$, $L_4=i_{511}$, while $ ~ n\ge 3$ is
a natural number, $|z|^2= zz^*$ and $|z| = \sqrt{zz^*}$,
$~|z|^2=\sum_jz_j^2$, $z=\sum_jz_ji_j$, $z_j\in {\bf R}$ for each
$j$, $z\in {\cal A}_v$. We take the branch of the square root
function $\sqrt{z}$ so that $\sqrt{b}>0$ for each $b>0$. The
conjugated number we write in the $z$-representation as $z^*=2 \pi
_0(z) -z$, where $\pi _0(z)=z_0$ is given by Formula $(3)$. That is,
the function $|z|$ is written in the $z$-representation. \par Since
the exponent series \par $(6)$ $e^z=\sum_{m=0}^{\infty }z^m/m!$
\\ converges for each $z\in {\cal A}_v$ and each function $\tau _l(z)$
is written in the $z$-representation, this function $g(z)$ is
analytic in the $z$-representation on the Cayley-Dickson algebra
${\cal A}_{v}$.
\par Each twice iterated exponent can be written in the form
$$(7)\quad \exp [ \tau _l(z) \exp (2\pi L_l \tau _l^{1/n}(z) ) ] =
\exp [ Re \{ \tau _l(z)\exp (2\pi L_l \tau _l^{1/n}(z) ) \} ] [ \cos
(|\tau _l(z)|)$$  $$ + \frac{ \sin (|\tau _l(z)|)} {|\tau _l(z)|} Im
\{ \tau _l(z) \exp (2\pi L_l \tau _l^{1/n}(z) ) \} ] =: e_l(z)$$ for
$\tau _l(z)\ne 0$, while $e^0=1$ for $\tau _l(z)=0$. Indeed, the
norm in the quaternion skew field ${\bf H}_l$ is multiplicative and
the quaternion skew field is without divisors of zero. Moreover,
$L_l\perp \tau _l^{1/n}(z)$ and the product $2\pi L_l \tau
_l^{1/n}(z)$ is purely imaginary and hence $|\exp (2\pi L_l \tau
_l^{1/n}(z) )|=1$, where $Im (z) =z- Re (z)$ is the imaginary part
of a Cayley-Dickson number $z$, $ ~ Re (z)=(z+z^*)/2$ is the real
part of $z$. Its value $e_l(z)$ belongs to the embedded copy ${\bf
H}_l$ of the quaternion skew field, ${\bf H}_l\hookrightarrow {\cal
A}_v$.
\par On the other hand, the intersection of three embedded copies of
the quaternion skew field ${\bf H}_1\cap {\bf H}_2\cap {\bf
H}_3={\bf R}$ is equal to the real field, which is the center
$Z({\cal A}_{v})={\bf R}$ of the Cayley-Dickson algebra ${\cal
A}_v$, also ${\bf H}_1\cap ({\bf H}_2{\bf H}_3)={\bf R}$. The last
term $| z - \sum_{l=1}^3 \upsilon _l(z) |^2 L_4$ in Formula $(5)$ is
orthogonal to other additives $\exp [ \tau _1(z) \exp (2\pi L_1 \tau
_1^{1/n}(z) ) ]$ and \\ $ - (\exp [ \tau _2(z) \exp (2\pi L_2 \tau
_2^{1/n}(z)) ]) (\exp [ \tau _3(z) \exp (2\pi L_3 \tau _3^{1/n}(z))
] )$ there. Therefore, $g(z)$ may be equal to zero only when all
three twice iterated exponents have real values.
\par Each number in the complex field $p \in {\bf C}_l$ is
orthogonal to the doubling generator $L_l$ of the quaternion skew
field ${\bf H}_l$ and in accordance with Formula $(7)$ the iterated
exponent $e_l(z)$ is real only when $|\tau _l^{1/n}(z)|\in {\bf Z}/2
:= \{ 0, \pm 1/2, \pm 1, \pm 3/2,... \} $, since $Re (L_l\tau
_l^{1/n}(z))=0$ and
\par $(8)$ $\exp (2\pi L_l \tau _l^{1/n}(z) ) = \cos (2\pi
|\tau _l^{1/n}(z)|) + \frac{ \sin (2\pi |\tau _l^{1/n}(z)|)} {|\tau
_l^{1/n}(z)|} L_l \tau _l^{1/n}(z)$. \\ The number $\tau _l(z)\in
{\bf C}_l$ is complex, consequently, $\tau _l^{1/n}(z)$ has $n$
distinct isolated roots in the complex field ${\bf C}_l$
corresponding to $n$ branches of the function $z^{1/n}$. Take for
definiteness the branch of $g(z)$ corresponding to the branch of the
$n$-th root such that $b^{1/n}>0$ for each $b>0$.
\par In accordance with the notation above
\par $(9)$ $g(z)=\exp [w\exp (2\pi L_1w^{1/n})]
- (\exp [x\exp (2\pi L_2x^{1/n})]) (\exp [y\exp (2\pi L_3y^{1/n})])
+ | z - \sum_{l=1}^3 \upsilon _l(z) |^2 L_4 =: \psi
(w,x,y;z-\sum_{l=1}^3 \upsilon _l(z))$, \\ where $w=\tau _1(z)$, $ ~
x=\tau _2(z)$ and $ ~ y = \tau _3(z)$.
\par This function $\psi (w,x,y;z-\sum_{l=1}^3 \upsilon _l(z))$
may have zeros only when all three complex variables $w_1:=w\in {\bf
C}_1$, $w_2:=x\in {\bf C}_2$ and $w_3:=y\in {\bf C}_3$ are real
$n$-th powers of half-integer or integer numbers $t_l\in {\bf Z}/2$,
$w_l=t_l^n$, while $z= \sum_{l=1}^3 \upsilon _l(z)$, since $\frac{
\sin (|t_l^n|)} {|t_l^n|} \ne 0$ is non-zero for each non-zero
$t_l\in ({\bf Z}/2)\setminus \{ 0 \} $, while $\lim_{\phi \to 0}
\frac{\sin (\phi )}{\phi }=1$ (see also Formulas $(7-9)$ above). It
can be lightly seen that $\psi (w,x,y;z-\sum_{l=1}^3 \upsilon
_l(z))=0$ for $w=0$ and $x = -y=t_2^n$ with $t_2\in {\bf Z}/2$ for
$n$ odd, also $x=0$ and $w=y=t_1^n$ with $t_1\in {\bf Z}/2$, also
$y=0$ and $w=x=t_1^n$ with $t_1\in {\bf Z}/2$, when $z= \sum_{l=1}^3
\upsilon _l(z)$. That is, the function $g(z)$ has only isolated
zeros in the Cayley-Dickson algebra ${\cal A}_{v}$. It is shown
below that this function has not any zeros when the product of
arguments $wxy\ne 0$ is non-zero for $n\ge 3$.
\par The equality $g(z)=0$ is equivalent to
\par $(10)$ $s(t_1) t_1^n = s(t_2) t_2^n + s(t_3) t_3^n$ and $z=
\sum_{l=1}^3 \upsilon _l(z)$, \\ where $t_1, t_2, t_3 \in {\bf
Z}/2$, $ ~ s(t_l)=1$ for an integer $t_l\in {\bf Z}$ and $s(t_l)=
-1$ for non-integer half-integer $t_l\in {\bf Z}+1/2$. If one of
these numbers is non-integer half-integer: either $t_1=k+1/2$ or
$t_2=l+1/2$ or $t_3=m+1/2$, where $k, l, m\in {\bf Z}$, then
Equality $(10)$ is equivalent to $s(t_1) (2k+1)^n = s(t_2) (2t_2)^n+
s(t_3) (2t_3)^n$ or $s(t_1) (2t_1)^n=s(t_2) (2l+1)^n+ s(t_3)
(2t_3)^n$ or $s(t_1) (2t_1)^n=s(t_2) (2t_2)^n+ s(t_3) (2m+1)^n$
respectively. That is Equality $(10)$ is always equivalent to the
corresponding equality for integers. \par If the number $n$ is even
and $(10)$ is satisfied, then it is also satisfied for $(\pm t_1,
\pm t_2, \pm t_3)$. If $n$ is odd and Equality $(10)$ is satisfied,
then it is also satisfied for the triple $(-t_1,-t_2,-t_3)$.
\par Considering different possible signs $s(t_1)$, $s(t_2)$, $s(t_3)$,
$sign (t_1)$, $sign (t_2)$ and $sign (t_3)$ one leads to the
conclusion that
$(10)$ is equivalent to the equality \par $(11)$ $a^n+b^n=c^n$ \\
for non-negative integers, where $c= \nu \max ( |t_1|, |t_2|,
|t_3|)$; $ ~ \nu =1$, when $t_1$ and $t_2$ and $t_3$ are integer; $
~ \nu =2$, when at least one of these numbers $t_1$ or $t_2$ or
$t_3$ is non-integer half-integer. Indeed, for non-negative numbers
$t_1$ and $t_2$ and $t_3$ the equality $t_1^n=t_2^n-t_3^n$ is
equivalent to $t_1^n+t_3^n=t_2^n$; $t_1^n= -t_2^n+t_3^n$ is
equivalent to $t_1^n+t_2^n=t_3^n$; the identity $-t_1^n =
t_2^n+t_3^n$ may be satisfied only for $t_1=t_2=t_3=0$; while
$-t_1^n=-t_2^n+t_3^n$ is equivalent to $t_1^n+t_3^n=t_2^n$; then
$-t_1^n=t_2^n-t_3^n$ is equivalent to $t_1^n+t_2^n=t_3^n$; the
equality $t_1^n = -t_2^n-t_3^n$ may be satisfied only for
$t_1=t_2=t_3=0$; also $-t_1^n=-t_2^n-t_3^n$ is equivalent to
$t_1^n=t_2^n+t_3^n$, where $t_1\ge 0$, $t_2\ge 0$ and $t_3\ge 0$.
Thus all zeros $(w,x,y)$ of the function $g(z)$ are described by the
set of all non-negative integer solutions of Equation $(11)$ with
the condition $z= \sum_{l=1}^3 \upsilon _l(z)$ up to the symmetry
transformations $(t_1,t_2,t_3)\mapsto (t_{\sigma (1)},t_{\sigma
(2)},t_{\sigma (3)} )$ and the multiplier $\nu \in \{ 1, 2 \} $ as
above, where $\sigma : \{ 1, 2, 3 \} \to \{ 1, 2, 3 \} $ is a
bijective surjective mapping.

\par Next it is verified that the residue operator of the reciprocal
function $f(z) := \frac{1}{g(z)}$ is non-degenerate for each zero
$(w, x, y)=(t_1^n, t_2^n, t_3^n)$, $~ t_1, t_2, t_3\in ({\bf
Z}/2\setminus \{ 0 \} ) $ of $g(z)$ if such zero exists. From the
proof above we know that all poles $\omega $ of the reciprocal
function $f(z)$ or equivalently zeros of the function $g(z)$ are
isolated points when zeros $z=\omega =(t_1^ni_{\kappa (1)}+
t_2^ni_{\kappa (2)} + t_3^ni_{\kappa (3)})$ with $t_1, t_2, t_3\in
({\bf Z}/2)$ exist.
\par The function $g(z)$
is analytic in the $z$-representation on the Cayley-Dickson algebra
${\cal A}_{v}$ and the reciprocal function $f(z)$ is meromorphic. In
this case residue operators of meromorphic functions of
Cayley-Dickson variables at point poles were defined and studied in
\cite{lufjmsrf,lujmsnfpcd,ludfov,ludoyst}.
\par Calculating the (super-)derivative operator of the function
$g(z)$ one obtains due to the  chain rule (see also Proposition 5
above):
\par $(12)$ $(dg(z)/dz).h = (d\exp [ \tau _1(z) \exp (2\pi
L_1 \tau _1^{1/n}(z) ) ]/dz).\upsilon _1(h) $\par $ - [(d(\exp [
\tau _2(z) \exp (2\pi L_2 \tau _2^{1/n}(z)) ])/dz).\upsilon _2(h)]
(\exp [ \tau _3(z) \exp (2\pi L_3 \tau _3^{1/n}(z)) ] )$\par $ -
(\exp [ \tau _2(z) \exp (2\pi L_2 \tau _2^{1/n}(z)) ]) [(d(\exp [
\tau _3(z) \exp (2\pi L_3 \tau _3^{1/n}(z)) ] )/dz).\upsilon _3(h)]
+ (d|z- \sum_{l=1}^3 \upsilon _l(z)|^2L_4/dz).(h- \sum_{l=1}^3
\upsilon _l(h)),$
\\ for each $z, h \in {\cal A}_{v}$, where
$$(13)\quad (d\exp [ \tau _l(z) \exp (2\pi L_l \tau _l^{1/n}(z) )
]/dz).\upsilon _1(h) =$$  $$ [\sum_{k=1}^{\infty } \sum_{m=0}^{k-1}
\{ \xi ^m I \xi ^{k-m-1} \} /k!  ]. (I e^{2\pi L_l\vartheta }$$
$$ + \vartheta ^n [\sum_{k=1}^{\infty } \sum_{m=0}^{k-1} \{ (2\pi
L_l \vartheta )^m (2\pi L_l I) (2\pi L_l \vartheta )^{k-1-m} \} /k!
]].[(d\vartheta /dz). \upsilon _l(h)] ,$$
\par $(14)$ $(d\vartheta /dz).h =
[d e^{(Ln \tau _l(z))/n}/dz].\upsilon _l(h)=\tau _l^{-1+1/n}(z)\tau
_l(h)/n$,
\par $(15)$ $(d|z-
\sum_{l=1}^3 \upsilon _l(z)|^mL_4/dz).(h- \sum_{l=1}^3 \upsilon
_l(h))= \frac{m}{2} \alpha ^{(m-2)/2} \{ (h- \sum_{l=1}^3 \upsilon
_l(h)) \{ 2 \pi _0(z) - z + \sum_{l=1}^3 \upsilon _l(z) \}  + (z -
\sum_{l=1}^3 \upsilon _l(z)) \{ 2\pi _0(h) - h +
\sum_{l=1}^3 \upsilon _l(h) \} \} L_4$ for each $m\in {\bf R}$, \\
since $\upsilon _l$ is the real-linear projection operator and $\tau
_l(z)=\tau _l(\upsilon _l(z))$ and $\pi _0(\upsilon _l(z))=0$ for
each $l=1, 2, 3$ and $z\in {\cal A}_{v}$, where each curled bracket
corresponds to the right order of multiplication $ \{ ab...cd \} =
a(b(...(cd)...))$, the letter $ ~I$ denotes the unit operator acting
here on $\mbox{}_lh := \upsilon _l(h)$, $ ~ \xi = \vartheta ^n
e^{2\pi L_l \vartheta }\in {\bf H}_l$, $ ~ \vartheta = \tau
_l^{1/n}(z)\in {\bf C}_l$, $ ~ |z|^m = \exp [\frac{m}{2} Ln (z(2\pi
_0(z)-z))]$, $ ~ \alpha = |z- \sum_{l=1}^3 \upsilon _l(z)|^2$ in
accordance with the chosen $z$-representation above with the help of
Formulas $(2,3)$ and 2$(1-3)$.

\par For definiteness we use the left algorithm for calculation of
line integrals of functions over Cayley-Dickson algebras. Then a
residue of a meromorphic function $q$ on a domain $U$ with a
singularity at an isolated point $\xi \in {\cal A}_r$ is defined as
$$(16)\quad Res_{\gamma _K} (\xi ,q).M :=
(2\pi )^{-1}\lim_{\delta \downarrow 0}(\int_{\gamma _{\delta
,K}}q(z)dz)$$ whenever this limit exists, where
\par $(17)$ $ \gamma _{\delta ,K}(t)=
(\xi + \delta \chi (t)K)\subset U\setminus \{ \xi \} \subset {\cal A}_r,$ \\
$\chi $ is a rectifiable closed path (i.e. loop) winding one time
around zero,
\par $(18)$ $2\pi M= \int_{\chi } dLn z$,
\\ $K$ is a marked Cayley-Dickson number for $\gamma $ such that
$|K|=1$, $~ 0< \delta $, $~ \gamma _K := \gamma _{1,K}$. Making the
change of the variable $z\mapsto zK$ one can relate the case of
$K=1$ with the case of $K\in {\cal A}_r\setminus {\bf R}$ in
Formulas $(16,17)$.
\par Particularly it may be the circle $\chi (t) = \rho \exp
(2\pi tM)$, where $\rho >0$, $M\in {\cal S}_r$, $t\in [0,1]$,
$0<\delta \le 1$, there is not any other singular point in the
closed ball $B({\cal A}_r, \xi ,\rho )$ of radius $\rho $ and the
center at $\xi $ in the Cayley-Dickson algebra ${\cal A}_r$, (see
also \S 1 in \cite{lufjmsrf}), here $r=v=9$, where ${\cal S}_r := \{
M\in {\cal A}_r: ~ Re (M)=0, ~ |M|=1 \} $ denotes the purely
imaginary unit sphere.
\par Let a purely imaginary Cayley-Dickson number be $M\in {\bf C}_1\cap
{\cal S}_r$ and $g(\xi )=0$, then from Formulas $(12-18)$ and
$\partial (\xi + \delta \rho \exp (2\pi tM)K) /\partial t = \delta
\rho 2\pi (M \exp (2\pi tM))K $ it follows that
$$(19)\quad Res_{\gamma _{i_{\kappa (1)}}} (\xi ,f).M =
[M e^{ - 2\pi L_1w^{1/n}} \exp (- w e^{2\pi L_1w^{1/n}})] [1 +
\frac{2\pi}{n} L_1 w^{1/n}]^{-1},$$ since $\alpha (\beta {\bf l}) =
(\beta \alpha ){\bf l}$ for each $\alpha , \beta \in {\cal A}_k$ and
${\bf l} = i_{2^k}$ for $k\ge 2$ in accordance with Formula 1$(2)$,
the quaternion skew field is associative and the octonion algebra
$alg _{\bf R}\{ M,L_1,i_{\kappa (1)} \} $ is alternative,
particularly for complex numbers $\alpha , \beta \in {\bf C}_M :=
{\bf R}\oplus {\bf R}M$, also $\upsilon _l\circ \upsilon _k(z)=0$
for each $k\ne l\in \{ 1, 2, 3 \} $, where $w=\tau _1(\xi )=t_1^n
\ne 0$, $x=\tau _2(\xi )=t_2^n$ and $y=\tau _3(\xi )=t_3^n\in {\bf
Z}/2$, the loop $\gamma _{i_{\kappa (1)}}$ is given by Formula
$(17)$. Symmetrically by substitution of variables for $M\in {\bf
C}_l\cap {\cal S}_r$ and non-zero $w_l\ne 0$ and $g(\xi )=0$ we get
$$(20) \quad Res _{\gamma _{i_{\kappa (l)}}} (\xi ,f).M = $$
$$- [M e^{ - 2\pi L_lw_l^{1/n}} \exp (- w_l e^{2\pi L_lw_l^{1/n}}) \exp
(- w_k e^{2\pi L_kw_k^{1/n}})] [1 +\frac{2\pi}{n} L_l
w_l^{1/n}]^{-1},$$ where $w_2=x$ and $w_3=y$, $(l,k)=(2,3)$ or
$(l,k)=(3,2)$, for $\gamma $ given by Formula $(17)$, $~ w=t_l^n$,
$t_l\in {\bf Z}/2\setminus \{ 0 \} $. Thus $z=wi_{\kappa
(1)}+xi_{\kappa (2)}+yi_{\kappa (3)} \in {\cal A}_v$ with $wxy\ne 0$
is a zero of the function $g(z)$ if and only if the residue operator
$Res _{\gamma i_{\kappa (l)}} (\xi ,f).M$ is non-degenerate by $M\in
{\bf C}_l\cap {\cal S}_r$.

\par If $M\in {\bf C}_1\cap {\cal S}_r$ and $n$ is even, one gets
from Formula $(19)$: $$(21)\quad Res_{\gamma _{i_{\kappa (1)}}} (\xi
,f).M + Res_{\gamma _{i_{\kappa (1)}}} (-\xi ,f).M =$$  $$  [2 M
e^{- 2 \pi L_1w^{1/n}} \exp ( -  w e^{2\pi L_1w^{1/n}}) ]/ [(4\pi ^2
w^{2/n}/n^2) +1],$$ if $n$ is odd:
$$(22)\quad Res_{\gamma _{i_{\kappa (1)}}} (\xi
,f).M + Res_{\gamma _{i_{\kappa (1)}}} (-\xi ,f).M = M [\exp (w
e^{2\pi L_1w^{1/n}}) + \exp (- w e^{2\pi L_1w^{1/n}})$$  $$+ \frac{2
\pi }{n} L_1w^{1/n} ( \exp (w e^{2\pi L_1w^{1/n}}) - \exp (- w
e^{2\pi L_1w^{1/n}}))] e^{2\pi L_1w^{1/n}}/ [(4\pi ^2 w^{2/n}/n^2)
+1].$$ From Formulas $(21)$ and $(22)$ we have, that  their sum by
$1/2\le t\le k$ for each $1/2\le k\in {\bf N}/2$ is non-zero, since
$e^{2\pi L_1w^{1/n}}\in \{ -1, 1 \} $, $ ~ t\in {\bf N}/2$, $ ~
w=t^n$. Analogously we get this conclusion for $Res_{\gamma
i_{\kappa (l)}}$ with $l=2, 3$ also with $[ - \exp (w_l e^{2\pi
L_lw_l^{1/n}}) \exp (w_k e^{2\pi L_kw_k^{1/n}})]$ instead of $\exp
(w e^{2\pi L_1w^{1/n}})$. \par Now we consider an analytic change of
the variable $z\in {\cal A}_v$ of special type: \par $(23)$ $\eta =
\eta(z) := (z-\sum_{l=1}^3
\upsilon _l(z))+ \sum_{l=1}^3 \tau _l(z)^{1/n} i_{\kappa (l)}$, \\
taking the branch of the $n$-th root function $z^{1/n}$ such that
$b^{1/n}>0$ for each $b>0$, where $\eta (z)$ is written in the
$z$-representation due to Formulas $(2,3)$ and 2$(1-3)$. Therefore,
one gets the equality:
\par $(23.1)$ $|\eta (z)|^2 = |z- \sum_{l=1}^3 \upsilon _l(z)|^2 +
\sum_{l=1}^3 | \tau _l(z)|^{2/n}$. \\
 The inverse transform of the Cayley-Dickson
variable is:
\par $(23.2)$ $z = z(\eta ) := (\eta - \sum_{l=1}^3
\upsilon _l(\eta )) + \sum_{l=1}^3 \tau _l(\eta )^n i_{\kappa (l)}$. \\
Therefore, the function $g(z)$ can be presented as the composition
of two (super-) differentiable over the Cayley-Dickson algebra
${\cal A}_v$ functions $g(z)=q(\eta (z))$, where
\par $(24)$ $q(\eta ) := \exp [ (\tau _1(\eta ))^n \exp (2\pi
L_1 \tau _1(\eta ) ) ]$\\ $ - (\exp [ (\tau _2(\eta ))^n \exp (2\pi
L_2 \tau _2(\eta )) ]) (\exp [ (\tau _3(\eta ))^n \exp (2\pi L_3
\tau _3(\eta )) ] ) + | \eta  - \sum_{l=1}^3 \upsilon _l(\eta )
|^2L_4 $, \\ $q(\eta )$ is written in the $\eta $ representation
with the help of Formulas $(2,3)$ and 2$(1-3)$. The only zeros of
$q(\eta )$ are $\eta \in {\cal A}_v$ satisfying two conditions: \par
$(i)$ $\eta =\sum_{l=1}^3 \upsilon _l(\eta )$ and \par $(ii)$ $\tau
_1^n(\eta ) = \tau _2^n(\eta )+\tau _3^n(\eta )$ with $\tau _l(\eta
)\in T^n_l\times ({\bf Z}/2)$, where $T^n_l$ denotes the set of $n$
roots $\sqrt[n]{1}$ of the unit in the complex field ${\bf C}_l$.
\par If $\alpha $ and $\beta $ are two Cayley-Dickson numbers
so that \par $(iii)$ $\alpha = \alpha _0 + \alpha _1M$ and $\beta =
\beta _0 + \beta _1 M$ with real coordinates $\alpha _0, \alpha _1,
\beta _0, \beta _1\in {\bf R}$ and a purely imaginary Cayley-Dickson
number $M\in {\cal S}_v$, $|\alpha |>0$, then \par $(iv)$ $Ln
(\alpha + \beta ) = Ln (\alpha ) + Ln (1 + \beta /\alpha )$, \\
since $\alpha $ and $\beta $ commute in such case.
\par Let $x=x_0+x_1M$, $y=y_1K$ and $z=z_0+z_1M\in {\cal A}_r$ be any Cayley-Dickson numbers
such that $x_0, x_1, y_1, z_0$ and $z_1\in {\bf R}$ are real; $ ~ M$
and $K\in {\cal S}_r$ are purely imaginary and orthogonal $Re
(MK^*)=0$. Suppose that closed rectifiable curves $x(t)$, $y(t)$ and
$z(t)$ are given, $t\in [0,2\pi ]$, such that \par $(v)$ $|z(t)|>
|x(t)+y(t)|$ for each $t\in [0,2\pi ]$, where \par $(vi)$
$z=|z|e^{Mt}$ with $t\in [0,2\pi ]$, i.e. $\cos (t) =z_0/|z|$,
$|z|^2=z_0^2+z_1^2$. Then we infer that \par $(vii)$ $x+y+z =
|x+y+z| e^{L\phi }$, \\ where $|z|^2 = [ (x_0+z_0)^2 + (x_1+z_1)^2 +
y_1^2]^{1/2}$, $ ~ \cos (\phi ) = (x_0+z_0)/|x+y+z|$, $ ~ L=L(t) =
[(x_1+z_1)M + y_1K]/(|x+y+z| \sin (\phi ) ) =  \pm ((x_1+z_1)M
+y_1K) / \sqrt{(x_1+z_1)^2+y_1^2}$. \par At first we consider the
particular case when \par $(vii.1)$ $y_1(t)=y_1$ is constant and $K$
and $M$ are constant and \par $(vii.2)$ the algebra $alg_{\bf
R}(M,K)$ over the real field generated by $M$ and $K$ is contained
in an alternative sub-algebra in ${\cal A}_r$.
\par Let $\eta (t)=\alpha _0(t)+\alpha _1(t)M+y_1K$ be a curve in ${\cal A}_r$ with real functions
$\alpha _0(t+\pi ) = - \alpha _0(t)$ and $\alpha _1(t+\pi ) = -
\alpha _1(t)$ for each $t\in [0,2\pi ]$ so that $\alpha
_0^2(t)+\alpha _1^2(t)>y_1^2$ for each $t$, for example, when
$\alpha _0(t)+\alpha _1(t)M$ is a circle with the center at zero of
radius $\rho >|y_1|$. Therefore, one gets $\eta (t+\pi ) = -\alpha
_0(t) - \alpha _1(t)M + y_1K = - (\alpha _0(t) + \alpha _1(t)M -
y_1K)$. Then we infer the equalities: \par $(viii.1)$ $\int_0^{2\pi
} d Ln ~ \eta (t) = \frac{1}{2} \int_0^{4\pi } d Ln ~ \eta (t) =
\frac{1}{2} [\int_0^{2\pi } d Ln ~ \eta (t) + \int_0^{2\pi } d Ln ~
\eta (t+\pi )] $\par $ = \frac{1}{2} [ \int_0^{2\pi }d \ln (\alpha
_0^2+\alpha _1^2+y_1^2) + \int_0^{2\pi } d \{ [\alpha _1M +y_1K]\phi
/\sqrt{\alpha _1^2 + y_1^2} \} $\par $ + \int_0^{2\pi } d \{ [\alpha
_1M - y_1K](\phi +\pi )/\sqrt{\alpha _1^2 + y_1^2} \} ] $ \par $=
\frac{1}{2} \int_0^{2\pi }d \ln (\alpha _0^2+\alpha _1^2+y_1^2) +
\int_0^{2\pi } d \{ \alpha _1M \phi /\sqrt{\alpha _1^2 + y_1^2} \} $
\\ using $(vii)$,
consequently, \par $(viii.2)$ $\int_0^{2\pi } d Ln ~ \eta (t)\in
{\bf R}\oplus M{\bf R}$.
\par On the other hand, we have the equality $(\alpha _0 + \alpha
_1M +y_1K) = (\alpha _0K^* + \alpha _1MK^* +y_1)K$, since $M$ and
$K$ are contained in the alternative sub-algebra in the
Cayley-Dickson algebra ${\cal A}_r$ by our supposition made above.
Then we have the equality $Ln ~ (\alpha _0(t)K^* + \alpha _1(t)MK^*
+y_1) = Ln ~ (\alpha _0(t) K^* + \alpha _1(t) MK^*) +  Ln ~
(1+y_1/(\alpha _0(t)K^* + \alpha _1(t)MK^*))$ for each $t\in [0,2\pi
]$,  consequently, the winding number around zero of the curve
$(\alpha _0(t) + \alpha _1(t)M +y_1K)$ is equal to that of $(\alpha
_0(t) + \alpha _1(t)M)$, since $|y_1|/|\alpha _0(t)K^* + \alpha
_1(t)MK^*|= |y_1|/|\alpha _0(t) + \alpha _1(t)M|<1$ for each $t\in
[0,2\pi ]$.
\par Now we consider the general case \par $(vii.3)$
$M=M(t)\in{\cal  S}_r$ and $K=K(t)\in {\cal S}_r$ and they are
orthogonal $K\perp M$, i.e. $Re (KM^*)=0$, for each $t\in [0,2\pi
]$.
\par Making transformations  and using Formulas $(vi,vii)$ we deduce that
\par $(ix)$ $Ln ~ (x(t)+y+z(t)) - Ln ~ z(t) = Ln \frac{|x(t)+y+z(t)|}{|z(t)|}
+L(t)\phi (t) - Mt =$ \par $\frac{1}{2} Ln [1+\frac{(z_0+x_0)^2 +
(x_1+z_1)^2+ y_1^2 - z_0^2-z_1^2}{z_0^2+z_1^2}] + M [
\frac{(x_1+z_1)}{\sqrt{(x_1+z_1)^2+y_1^2}} \arccos
\frac{x_0+z_0}{|x+y+z|} - \arccos \frac{z_0}{|z|}] + K
\frac{y_1}{\sqrt{(x_1+z_1)^2+y_1^2}}\arccos \frac{x_0+z_0}{|x+y+z|}
$, \\ where $x_0=x_0(t)$, $x_1=x_1(t)$, $z_0=z_0(t)$, $z_1=z_1(t)$.
\par In view of Rouch\'e's theorem \cite{lavrsch,ludfov,ludoyst} in the complex plane
${\bf R}\oplus {\bf R}M$
\par $(x)$ $\int_0^{2\pi } d [Ln ~ (q(t)+z(t)) - Ln ~ z(t)) =
\int_0^{2\pi } d \{ \frac{1}{2} Ln [1+\frac{(z_0+q_0)^2 +
(q_1+z_1)^2 - z_0^2-z_1^2}{z_0^2+z_1^2}] + M [ \arccos
\frac{q_0+z_0}{|q+z|} - \arccos \frac{z_0}{|z|} ] \} =0$, \\ where
$q=q(t)=q_0(t)+q_1(t)M$ is a rectifiable closed curve (i.e. a loop:
$q_0(0)=q_0(2\pi )$) such that $|q(t)|<|z(t)|$ for each $t\in
[0,2\pi ]$. \par Conditions $(vii.1)$ and $(vii.2)$ can be abandoned
with the help of the homotopy theorem (see Theorem 2.15 in
\cite{ludfov,ludoyst}). It can be lightly seen that $(viii.1,x)$
imply the equality
\par $(xi)$ $\int_0^{2\pi } d \{ \frac{1}{2} Ln [1+\frac{(z_0+x_0)^2
+ (x_1+z_1)^2+ y_1^2 - z_0^2-z_1^2}{z_0^2+z_1^2}] $\par $ + M [
\frac{(x_1+z_1)}{\sqrt{(x_1+z_1)^2+y_1^2}} \arccos
\frac{x_0+z_0}{|x+y+z|} - \arccos \frac{z_0}{|z|}]$\par $ + K
\frac{y_1}{\sqrt{(x_1+z_1)^2+y_1^2}}\arccos
\frac{x_0+z_0}{|x+y+z|} \} =0$, \\
since $ |\frac{(z_0+x_0)^2 + (x_1+z_1)^2+ y_1^2 -
z_0^2-z_1^2}{z_0^2+z_1^2}|< 1$ for each $t$ and \par $\int_0^{2\pi }
d Ln [1+\frac{(z_0+x_0)^2 + (x_1+z_1)^2+ y_1^2 -
z_0^2-z_1^2}{z_0^2+z_1^2}]  =0$.
\par Indeed, consider the equation:
\par $\xi \cos (\alpha ) - \sqrt{1-\xi ^2} \sin (\alpha ) = \eta $,
\\ where $\alpha \in {\bf R}$, $~ \eta \in {\bf R}$, $~ |\eta |\le
1$. This gives $(\xi \cos (\alpha ) - \eta )^2 = (1-\xi ^2)\sin
^2(\alpha )$ or $\xi ^2 - 2\xi \eta \cos (\alpha ) + (\eta ^2 - \sin
^2(\alpha )) =0$. By Vieta's formula one gets two solutions $\xi
_{1,2} = \eta \cos (\alpha ) \pm \sqrt{1-\eta ^2} \sin (\alpha )$.
Particularly we take $\eta = \frac{1}{2\pi } \arccos
(\frac{z_0}{|z|})$, $~ \cos (\alpha ) =
\frac{x_1+z_1}{\sqrt{(x_1+z_1)^2+y_1^2}},$ $~ \sin (\alpha )
=\frac{y_1}{\sqrt{(x_1+z_1)^2+y_1^2}}.$ Therefore, we deduce that
\par $M[\frac{x_1+z_1}{\sqrt{(x_1+z_1)^2+y_1^2}} \arccos
(\frac{x_0+z_0}{|x+y+z|}) - \arccos (\frac{z_0}{|z|})] + K
\frac{y_1}{\sqrt{(x_1+z_1)^2+y_1^2}}\arccos
(\frac{x_0+z_0}{|x+y+z|})$\par $ = M [\arccos
(\frac{x_0+z_0}{|x+y+z|}) - \arccos (\frac{z_0}{|z|})]\cos (\alpha )
 + [K \arccos
(\frac{x_0+z_0}{|x+y+z|}) - 2\pi M \sqrt{1-\eta ^2}] \sin (\alpha
)$. \par The transformation ${\xi \choose {\sqrt{1-\xi ^2}}}\mapsto
{{\cos (\alpha ) - \sin (\alpha )}\choose{ \sin (\alpha ) ~ \cos
(\alpha ) }} {\xi \choose {\sqrt{1-\xi ^2}}}$ is orthogonal. In
Formula $(x)$ one can choose $\frac{x_0+z_0}{|x+y+z|} =
\frac{q_0+z_0}{|q+z|}$ with $|q|=|x+y|$ and $Im (q)\in Im(z){\bf
R}$, where $Im (q)=q-Re (q)$, $ ~ |x+y|^2=|x|^2+|y_1|^2$, since $x$
and $y$ are orthogonal, $x\perp y$. If $y_1=0$, then $|\cos (\alpha
)|=1$ and $\sin (\alpha )=0$. If $|y_1|>0$, then $|\cos (\alpha
)|<1$ and $|\sin (\alpha )|>0$. Using the homotopy theorem and
Formulas $(viii.1,viii.2, x)$ we infer that
\par $(xi.1)$ $\int_0^{2\pi } d \{ M [\arccos
(\frac{x_0+z_0}{|x+y+z|}) - \arccos (\frac{z_0}{|z|})]\cos (\alpha
)$ \par $ + [K \arccos (\frac{x_0+z_0}{|x+y+z|}) -  M \sqrt{(2\pi
)^2- \arccos ^2 (\frac{z_0}{|z|})}] \sin (\alpha ) \} (t) =0 $, \\
since the function in the integral does not change its branch,
consequently, the function $[Ln ~ (x(t)+y(t)+z(t)) - Ln ~ z(t)]$
does not change its branch.
\par Thus Formula $(xi.1)$ implies that
\par $(xii)$ ${\hat I}n (0,x+y+z) = {\hat I}n (0,z) := \frac{1}{2\pi }\int_0^{2\pi } d Ln
~ z(t)$, \\ when Conditions $(v,vi,vii.3)$ are satisfied (see also
Corollary 3.26.3 in \cite{ludfov,ludoyst}), where $\frac{1}{2\pi
}\int_0^{2\pi } d Ln ~ z(t)=P\in {\bf Z}{\cal S}_r$, particularly
$P=M$ for a constant $M\in {\cal S}_r$ (see above).

\par If $u=(u_1,...,u_k)$ is a vector in the Euclidean space ${\bf R}^k$,
then $|u|_p := \sqrt[p]{u_1^p+...+u_k^p}$ for $1\le p<\infty $. Let
$t>s>0$ and $0<a\le b$, then
\par $(25)$ $(a^t+b^t)^{1/t}< (a^s + b^s)^{1/s}$. \\
Indeed, Inequality $(25)$ is equivalent to $(1+x^t)^{1/t}
<(1+x^s)^{1/s}$ for each $x\ge 1$ or to
\par $(25.1)$ $(1+y^p)^{1/p} <1+y$, \\ where $x=b/a$, $ ~ p=t/s$, $ ~ y=x^s\ge
1$. Consider the function $g(y) = (1+y^p)^{1/p} - (1+y)$. One has
$g(1)=2^{1/p} - 2<0$ and $g'(y) = (1+y^p)^{-1 +1/p}y^{p-1} -1<0$ for
each $y\ge 1$, since $y^p <1+y^p$ $ ~ \Leftrightarrow $ $ ~ y^{p-1}
< (1+y^p)^{(p-1)/p}$ $\Leftrightarrow $ $(1+y^p)^{-1+1/p}
y^{p-1}-1<0$. Therefore, $g(y)<0$ for each $y\ge 1$, which implies
Inequality $(25)$.
\par Particularly,
\par $(26)$ $c^2<a^2+b^2$  and \par $(27)$ $c^{2n} > a^{2n}
+b^{2n}$ for any non-zero real numbers $a, b, c\ne 0$ such that
$|c|^n=|a|^n+|b|^n$ with $n\ge 3$.
\par Since $2^n > 1^n +1^n$ and $3^n>2^{n+1}$ for $n\ge 3$,
we consider the equation $a^n+b^n=c^n$ for half-integer numbers $a,
b, c\in ({\bf Z}/2) \setminus \{ 0 \} $ with $\max (|a|, |b|, |c|)
\ge 3$. That is in the ball $B({\cal A}_v,0,2^{n+1})$ in the
Cayley-Dickson algebra Equation $(10)$ with non-zero half-integer
numbers $t_1, t_2, t_3$ has not any solution, i.e. the function
$g(z)$ is non-zero, $g(z)\ne 0$, for the argument $z=t_1^n i_{\kappa
(1)}+ t_2^ni_{\kappa (2)} + t_3^ni_{\kappa (3)}\in B({\cal
A}_v,0,2^{n+1})$ with $t_1t_2t_3\ne 0$,  since $|z|^2 = t_1^{2n}+
t_2^{2n}+ t_3^{2n}$ for such $z$, where $B({\cal A}_v,\xi ,\rho ) :=
\{ z\in {\cal A}_v: ~ |z-\xi |\le \rho \} $, $~ \xi \in {\cal A}_v$,
$ ~ \rho >0$. \par We proceed by induction. Suppose that in the ball
$B({\cal A}_v,0, 2 \rho ^n)$ the function $g(z)$ is no-zero for
$z=c^n i_{\kappa (1)}+a^ni_{\kappa (2)} + b^ni_{\kappa (3)}\in
B({\cal A}_v,0,2\rho ^n)$ with $abc\ne 0$, $ ~ a, b , c \in ({\bf
Z}/2)\setminus \{ 0 \} $, where $\rho \in {\bf N} := \{ 1, 2, 3,...
\} $ is a natural number.
\par We take the path $\kappa (t)$ consisting of two circles $\gamma ^1$
and $\gamma ^2$ of radius $(\rho  + 1 + \epsilon )^n$ and $(\rho  +
\epsilon )^n$, where $0<\epsilon <1/4$, and a joining them path
$\omega $ gone twice in one and the opposite direction such that
$\kappa (t)$ contains no any Cayley-Dickson number $z$ with
half-integer coordinates $z_j$ for any $j=0,...,2^v-1$. That is,
$\kappa (t)=\gamma ^1(4t)$ for $0\le t<1/4$, $\kappa (t)=\omega (
4t-1)$ for $1/4\le t<1/2$, $\omega (0)=\gamma ^1(1)$, $\kappa
(t)=\gamma ^2(3-4t)$ for $1/2\le t <3/4$, $\kappa (t)=\omega (4-4t)$
for $3/4\le t\le 1$, $\omega (1)=\gamma ^2(0)$, where $ \gamma ^1(t)
= (\rho +1 + \epsilon )^n \chi (t)i_{\kappa (1)},$ $ ~ \gamma ^2(t)
= (\rho  + \epsilon )^n \chi (t)i_{\kappa (1)},$ $ ~ \chi (t) = \exp
(2\pi tM)$, where $\rho >0$, $M=i_1$, $t\in [0,1]$.
\par In view of the theorem about change of a variable in the line
integral over the Cayley-Dickson algebra ${\cal A}_v$ and the
theorem about residues (see Theorems 2.6, 2.11 and 2.13 in
\cite{lufjmsrf} and \cite{ludfov,ludoyst})  we infer that
$$(28)\quad \int_{\kappa }f(z)dz = \int_{\psi } \frac{1}{q(\eta
)}d\eta ,$$ where $\psi (t) = \eta (\kappa (t))$ for each $t\in
[0,1]$, consequently,
$$(29)\quad \sum_{\xi =
w_1^ni_{\kappa (1)}+w_2^ni_{\kappa (2)}+w_3^ni_{\kappa (3)}; ~ (\rho
+\epsilon )\le |c|<(\rho +1+\epsilon ); ~ c^n=a^n+b^n} Res_{\gamma
_{i_{\kappa (1)}}} (\xi ,f).M $$
$$ = \sum_{\xi =
w_1^ni_{\kappa (1)}+w_2^ni_{\kappa (2)}+w_3^ni_{\kappa (3)}; ~ (\rho
+1+\epsilon )^2 > w_l^2 + w_k^2; ~ |c| \ge (\rho +\epsilon ); ~
c^n=a^n+b^n} Res_{\gamma _{i_{\kappa (1)}}} (\xi ,f).M$$ due to
Formulas $(23-24)$ and $(i,ii,xii)$ above, where $|c| = \max
(|w_1|,|w_2,|w_3|) = |w_m|$, these numbers are ordered as $|w_l|\le
|w_k|<|c|$, $ ~ c = w_m$, $ ~ a=w_l$, $ ~ b=w_k$, $~a, b, c\in {\bf
Z}/2$, $~l, k, m$ are pairwise distinct, $~ l\ne k \ne m$, $~ l\ne
m$; $ ~ l, k$ and $m \in \{ 1, 2, 3 \}$, since the algebra $alg_{\bf
R} (M,i_{\kappa (2)},i_{\kappa (3)})$ is alternative.
\par But Formulas $(26-29)$ give us the contradiction with the
supposition that there are half-integer solutions $c^n=a^n+b^n$ for
$n\ge 3$ with the non-zero product $abc\ne 0$, since $0<\epsilon
<1/4$ may be arbitrary small and the limits $\lim_{\downarrow
\epsilon }$ of the right and the left sides of $(29)$ coincide. Thus
in the ball $B({\cal A}_v,0,2(\rho +1)^n)$ there is not any solution
of the equation $g(\xi )=0$ with $\xi = w_1^ni_{\kappa
(1)}+w_2^ni_{\kappa (2)}+w_3^ni_{\kappa (3)}$ and $w_1, w_2, w_3\in
({\bf Z}/2)\setminus \{ 0 \} $, since $|\xi |^2 =
a^{2n}+b^{2n}+c^{2n}=w_1^{2n}+w_2^{2n}+w_3^{2n}$. The proof above
leads to the conclusion that Equation $(11)$ has not any solution in
natural numbers $a, b, c \in {\bf N}$ for $n\ge 3$, since $\rho \ge
2$ is arbitrary.
\par {\bf 8. Remark.} Alternatively it is possible to use
the argument principle and the Cayley-Dickson analogs of Rouch\'e's
theorem and the homotopy theorem instead of residues in \S 7. For
this purpose additional functions can be used.
\par  Denote by $\theta _{k,l} :{\bf H}_k\to {\bf H}_l$ isomorphisms
of copies of the quaternion skew field. Making the substitution of
variables and using the isomorphisms $\theta _{k,l}$ and the
equation $w^n=x^n+y^n$, we write two new functions $g_1$ and $g_2$
which have the same zeros as $g$:
\par $(1)$ $g_1(z) = g(z-\upsilon _2(z)-\upsilon _3(z) +
[(\theta _{1,2}(\tau _1(z)) - (\theta _{3,2}(\tau _3(z))] i_{\kappa
(2)} +[(\theta _{1,3}(\tau _1(z)) - (\theta _{2,3}(\tau _2(z))]
i_{\kappa (3)})$ and
\par $(2)$ $g_2(z) = g(z-\upsilon _1(z) -\upsilon _2(z) +
[(\theta _{2,1}(\tau _2(z)) + (\theta _{3,1}(\tau _3(z))]i_{\kappa
(1)}+ [(\theta _{1,2}(\tau _1(z)) - (\theta _{3,2}(\tau _3(z))]
i_{\kappa (2)})$ and
\par $(3)$ $g_3(z) = g(z-\upsilon _1(z) -\upsilon _3(z) +
[(\theta _{2,1}(\tau _2(z)) + (\theta _{3,1}(\tau _3(z))]i_{\kappa
(1)}+ [(\theta _{1,3}(\tau _1(z)) - (\theta _{2,3}(\tau _2(z))]
i_{\kappa (3)})$ \\
and analogous transformed functions $q_l(\eta )$, $l=1, 2, 3$.

\par Department of Applied Mathematics,
\par Moscow State Technical University MIREA,
av. Vernadsky 78,
\par Moscow, Russia
\par e-mail: sludkowski@mail.ru
\end{document}